\documentstyle[12pt]{amsart}
\title[Representation theory of Lie algebras]{Geometric 
representation theory of restricted Lie
algebras of classical type}
\author{Ivan Mirkovi\'{c}}
\address{Dept. of Math. and Stat., LGRT, UMass, Amherst, 
MA, 01003, USA}
\email{mirkovic@@math.umass.edu}
\author{Dmitriy Rumynin}
\address{
Mathematics Dept., University of Warwick,
Coventry, CV4 7AL, U.K.}
\email{rumynin@@maths.warwick.ac.uk}
\date{July 29, 1998; revised on June 20, 1999}
\thanks{The research was supported by NSF
and completed at University of Massachusetts at Amherst}

\newcommand{\r}{\rho_\eta}

\newcommand{\mo}{{\mathcal O}}
\newcommand{\ii}{\mbox{Id}}

\newcommand{\KK}{{\Bbb K}}

\newcommand{\BF}{{\mathbb F}}
\newcommand{\BL}{{\mathbb L}}
\newcommand{\BG}{{\mathbb G}}
\newcommand{\B}{{\mathcal B}}

\newcommand{\fl}{{\frak l}}

\newcommand{\fA}{{\mathfrak A}}
\newcommand{\fP}{{\mathfrak P}}

\newcommand{\ug}{U_{\chi}({\frak g})}
\newcommand{\ul}{U_{\chi}({\frak l})}

\newcommand{\spr}{{\mathcal B}^{\chi}}
\newcommand{\ver}{Z_{\chi,{\mathfrak b}}(\lambda)}

\newcommand{\g}{{\frak g}}

\newcommand{\fb}{{\mathfrak b}}

\newcommand{\BB}{{\mathcal B}}
\newcommand{\mA}{{\mathcal A}}
\newcommand{\LL}{{\mathcal L}}

\newcommand{\OO}{{\mathcal O}}
\newcommand{\GG}{{\mathcal G}}
\newcommand{\FF}{{\mathcal F}}

\newcommand{\fL}{{\frak L}}
\newcommand{\Ni}{\mbox{Ni}_\eta(c)}

\begin{document}

\begin{abstract}
We modify the Hochschild $\varphi$-map to construct central
extensions of a restricted Lie algebra. 
Such central extension
gives rise to a group scheme which leads to
a geometric construction of 
unrestricted representations. 
For a  classical semisimple
Lie algebra,  
we construct equivariant line bundles
 whose global sections afford
representations
with a nilpotent $p$-character.
\end{abstract}

\maketitle

Let $G$ be a connected simply connected semisimple algebraic group
over an algebraically closed field $\KK$ of characteristic $p$ and 
$\g$ be its Lie algebra.
The representation theory of $\g$ is connected with
the coadjoint orbits through
the notion of
a
$p$-character \cite{kac1,fri,jan4,hum4}.
An irreducible representation $\rho$ is finite-dimensional
and 
determines a $p$-character $\chi \in \g^{\ast}$ by
$\chi (x)^p\,\mbox{Id} = \rho (x)^p- \rho (x^{[p]})$
for each $x \in \g$ \cite{kac1}. 
There are indications
that a geometry stands behind 
this representation theory,
for instance, 
the Kac-Weisfeiler conjecture 
proved by 
Premet \cite{pre}.  
This work has been motivated by 
an idea of Humphreys 
that 
the representations affording $\chi$
should be related to 
the Springer fiber $\spr$.
Some of our  intuition comes from algebraic calculations of
Jantzen \cite{jan2,jan3}.
The most interesting evidence
for the relation between
Springer fibers and representations
of $\g$ 
is now given by
Lusztig  \cite{lus3}.

The main goal of 
this paper is to introduce a method for 
constructing unrestricted representations 
of $\g$ by taking global
sections of line bundles on infinitesimal neighborhoods
of certain subvarieties of $\spr$.
A more general approach implementing twisted sheaves of 
crystalline differential operators will be explained
elsewhere.

An attempt to study representations of $\g$ with a single
$p$-character $\chi$ has led to the notion of a reduced
enveloping algebra.
We modify this approach by considering a set
of $p$ different $p$-characters
$\{0,\chi,2\chi, \,\ldots\, ,
\newline
(p-1)\chi\}$ 
together in Section 1.
The category of such representations
is closed under  tensor products.
These are  restricted representations
of a central extension $\g_\chi$ of $\g$ by 
the multiplicative restricted Lie algebra $\g_m$.
One can think of this construction as a multiplicative version
of the Hochschild $\varphi$-map.

We discuss a geometric machinery necessary for the construction
of representations in Sections 2 and 3. In Section 4,  
we introduce equivariant line bundles
and construct representations. 
This section contains the main result 
(Theorem 4.3.2) of this paper, which
is a geometric construction of unrestricted representations.
Section 5 is devoted to various comments 
on the representations constructed.

Let us briefly explain the construction.
The central extension $\g_\chi$ defines a central extension
$
0
\rightarrow 
G_m^1
\rightarrow 
G^\chi
\rightarrow 
G^1
\rightarrow 
0
$
of the Frobenius kernels of
$G$ and the multiplicative group
$G_m$.
The group scheme $G^\chi$ acts on the flag variety $\B$
and preserves the
Frobenius neighborhood $\widehat{Z}$ of 
any subscheme $Z$.
For a $G$-equivariant line bundle $\FF_\lambda$
on $\B$,   
we 
construct 
a $G^\chi$-action
on 
$\FF_\lambda|_{\widehat{Z}}
$ 
with a 
``central charge $1$''. Then
 $\g$ will act on the global sections of
$\FF_\lambda|_{\widehat{Z}}
$ 
with a $p$-character $\chi$.

It suffices to construct such
an equivariant structure on a subscheme $\widehat{X}$ 
that contains
$Z$. 
We want to choose $X$ so that we can put hands on
the Frobenius neighborhood $\widehat{X}$.
We will assume
that $X$ is smooth 
so that
$G^\chi \times X\rightarrow \widehat{X}$ is
the quotient
map by the
action groupoid $\GG^\chi_X$
arising from the
$G^\chi$-action
on $\widehat X$.

To construct an equivariant structure, 
it suffices to split
the groupoid 
$\GG^\chi_X$ 
as
a product of the Frobenius kernel of the multiplicative
group $G_m^1$ and another groupoid $\GG^1_X$.
A necessary condition for this construction
is that
$X$ is a subvariety of $\spr$.

The groupoid
$\GG^\chi_X$ 
splits canonically over the diagonal.
We linearize the requirement
that this
splitting 
extends off the diagonal,
and study it in terms of Lie
algebroids of the above mentioned groupoids. 

The authors are greatly indebted 
to J. Humphreys whose inspiration
was crucial for writing this article. 
The authors would like to thank 
T. Ekedahl, J. Jantzen, J. Paradowski,  
G. Seligman, and  S. Siegel  
for various information.

\subsection{Notational conventions}
Let $\BF$ be the prime subfield
of an algebraically closed field $\KK$ 
of characteristic $p$.

\subsubsection{Restricted Lie algebras}
The main object of our study is
a finite dimensional restricted Lie algebra $\fl$
over $\KK$. If 
$\fl$ is the Lie algebra of a linear algebraic
group, the group is denoted by $L$.
While discussing a semisimple algebraic group,
we denote the group by $G$ and its Lie algebra by $\g$. 
Let $R_\g$ be the set
of roots of $\g$, $\Delta_\g$ 
be a set of simple roots, 
$W$ be the Weyl group, 
and $\Pi$ be the weight lattice.
The multiplicative
group and its Lie algebra are denoted 
$G_m$ and $\g_m$.

\subsubsection{Flag variety}

Let $\BB$ be the flag 
variety of $G$. We think of points
of $\BB$ over $\KK$ 
as 
Borel
subalgebras $\fb$ in $\g$.
Let ${\mathcal B}_w$ be a Schubert variety
for $w\in W$.
If $\chi$ is nilpotent then 
{\em the Springer fiber}
$ {\mathcal B}^{\chi}$ is a reduced subscheme of $\BB$,
whose points over $\KK$ are those Borel
subalgebras on which $\chi$ vanishes.

\subsubsection{Enveloping algebras}
The universal enveloping algebra of $\fl$ is
$U(\fl)$.
It contains a central Hopf
subalgebra $O$ generated by 
$x^p-x^{[p]}$ for all $x \in \fl$.
For any $\chi \in \fl^\ast$,  
the reduced enveloping algebra $\ul$
is a quotient of $U(\fl)$ by 
the ideal generated by $x^p-x^{[p]}- 
\chi (x)^p1_{U(\fl)}$ for all $x \in \fl$ \cite{far,fri}.  
The reduced enveloping algebra $U_0(\fl)$ is 
the restricted enveloping algebra $u(\fl)$.
All $\ul$ 
are twisted products of $u(\fl)$ with the field $\KK$
\cite{rum}.

\subsubsection{$p$-character}
A representation of $\fl$ 
has a $p$-character $\chi\in\fl^\ast$ if the representation
determines a $\ul$-module.
While working with $\g$, we 
assume that  
$\chi$ is a nilpotent element of $\g^\ast$. 
The case of a general $\chi\in\g^\ast$
can be reduced to the nilpotent case.

\subsubsection{Induction}

If ${\mathfrak m}$ is a restricted 
Lie subalgebra of $\fl$ such that
$\chi |_{\mathfrak m}=0$
then 
the induction functor 
Ind$_{U_\chi({\mathfrak m})}^{U_\chi(\fl)}$
is defined on the category of left 
$u({\mathfrak m})$-modules.
 
In  particular,
for a Borel subalgebra $\fb$ to $\g$,
if $\chi |_\fb=0$
 then all simple 
modules over $U_\chi (\fb)=u(\fb)$ 
are one-dimensional and parametrized
by the reduced (modulo $p$) weight lattice $\Lambda$.
The induced module
$
\ver = \ug \otimes_{u(\fb)} {\KK}_{\lambda}
$, $\lambda \in \Lambda$, 
called {\em a baby Verma module}, 
was introduced in \cite{hum0}.
Any irreducible $\ug$-module 
is a quotient of at least one  $\ver$, though
the module $\ver$ need not 
have a unique simple quotient, which
makes a classification of simple $\g$-modules
an interesting problem
\cite{fri}. 

\section{Central extensions}

\subsection{Central extensions of Hopf algebras}
Our approach will be explained 
in this section. The ground
field {\bf k} is arbitrary for this section. 

\subsubsection{}
Let us consider a Hopf algebra 
$U$ and its central Hopf subalgebra $O$.
Given $\chi\in\mbox{Spec}\, O({\bf k})$,
representations  in which $O$
acts by  $\chi$ 
are those that can be reduced 
to the algebra $U_{\chi}=U \otimes_O
{\bf k}({\chi})$. The algebra $U_{\chi}$ 
is not necessarily a Hopf algebra.
The basic idea of 
the present paper is to replace the study 
of $U_{\chi}$-modules
for a single $\chi$ with the study of 
$U_{\chi}$-modules as $\chi$ runs
over a closed subgroup of $\mbox{Spec} \, O$. 
One has more modules
but we benefit from having 
a Hopf algebra rather than 
a Hopf-Galois extension.

\subsubsection{Proposition} \label{hopf}
Let $O \rightarrow R$ be the natural 
map where $R$ is the algebra of functions
on the closed subgroup scheme 
of Spec$\, O$ generated 
by $\chi$. Then $U \otimes_O R$
is a Hopf algebra.
 
\subsubsection{Proof}
A subgroup scheme $X$
gives rise to a surjective 
Hopf algebra map
$\pi:O \rightarrow \mo(X)$.
We need to check that $A \otimes_O \mo(X)$
is a Hopf algebra.

The tensor product 
$C=A \otimes_{\bf k} \mo(X)$ is obviously a Hopf algebra. 
It suffices to check that
the ideal $I$, generated by all 
$x \otimes 1-1 \otimes \pi(x)$ 
with $x \in O$, is a Hopf ideal.
The latter means that the quotient
$C/I=A \otimes_O \mo(X)$
admits a Hopf algebra 
structure such that the quotient map
is a Hopf algebra homomorphism. 
Being a Hopf ideal includes three
axioms that we are checking now.

{\em Axiom 1:} $\varepsilon (I)=0$.
A typical element of $I$ has a form
$\sum_i a_i(x_i \otimes 1-1 \otimes \pi(x_i))b_i$
where $a_i,b_i\in C$, $x_i\in O$. Now we compute
$\varepsilon(\sum_i a_i
(x_i \otimes 1-1 \otimes \pi(x_i))b_i) =
\sum_i \varepsilon(a_i)(\varepsilon(x_i) \otimes 1-1 
\otimes \varepsilon(\pi(x_i)))\varepsilon(b_i) = 
\sum_i \varepsilon(a_i)(\varepsilon(x_i)1_C - 
\varepsilon(x_i)1_C)\varepsilon(b_i) = 0.$

{\em Axiom 2:} $S(I)\subseteq I$. Let us compute
$S(\sum_i a_i(x_i \otimes 1-1 \otimes \pi(x_i))b_i) =
\sum_i S(b_i)(S(x_i) \otimes 1-1 
\otimes S(\pi(x_i)))S(a_i) = 
\sum_i S(b_i)(S(x_i) \otimes 1-1 
\otimes \pi(S(x_i)))S(a_i) \in I.$

{\em Axiom 3:} $\Delta(I)\subseteq C\otimes I+I\otimes C$. 
Let us compute
$\Delta(\sum_i a_i(x_i \otimes 1-1 \otimes \pi(x_i))b_i) =
\sum_i \Delta(a_i)(x_{i1} \otimes 1 
\bigotimes x_{i2}\otimes 1
-1 \otimes \pi(x_{i1}) \bigotimes 
1\otimes\pi(x_{i2}))\Delta(b_i)
=
\sum_i \Delta(a_i)
[\{x_{i1} \otimes 1-1\otimes\pi(x_{i1})\}
\bigotimes x_{i2}\otimes 1
+
1\otimes\pi(x_{i1}) \bigotimes 
\{x_{i2}\otimes 1-1 \otimes \pi(x_{i2})\}]\Delta(b_i)
\in I\otimes C+C\otimes I.$
 $\ \ \ \ \Box$

\subsubsection{}
In the present paper, 
we focus on the case 
of the universal enveloping Hopf algebra
of a restricted Lie algebra. 
The subgroup generated by $\chi$ is ${\Bbb F}\chi$. 
A quantum linear group $O_q(G)$ and the unrestricted form 
of a quantum enveloping 
algebra $U_q(\g)$ at a root of unity
are other interesting options \cite{rum}.
However, 
a closed subgroup of $G$ or 
${\Bbb C}^{n-r} \times {\Bbb C}^{\ast r}$
generated by an element is
more complicated.

\subsection{Extensions of restricted Lie algebras.}

\subsubsection{} 
An exact sequence of restricted Lie algebras 
$0\rightarrow
{\frak a} \rightarrow 
{\frak m} \rightarrow \fl \rightarrow 0$ 
is called
{\em a central extension} of $\fl$ if $\frak a$ 
is a central ideal of ${\frak m}$. 
This terminology is not standard.
An additional condition ${\frak a}^{[p]}=0$ 
is required in \cite{ho2}
for a central extension.
The reason for this constraint is that such central
extensions can be parametrized 
by the second restricted cohomology group.
The important choice of $\frak a$ for us is $\g_m$,
which means that we usually have
${\frak a}^{[p]}={\frak a}$.

\subsubsection{Multiplicative Hochschild $\varphi$-map}

The Hochschild $\varphi$-map \cite{ho1,sul} provides
a central extension of $\fl$ by the additive Lie algebra
for each  $\chi \in \fl^{\ast}$.
We modify this construction to obtain a central extension
by $\g_m$ instead.
Given $\chi \in \fl^{\ast}$, we construct a central extension
$\fl_\chi$. This extension 
is trivial as an extension of Lie algebras,
i.e.
$\fl_{\chi}=\fl \oplus {\KK}c$;
but the $p$-structure is twisted by $\chi$:
\begin{equation}
(a+ \alpha c)^{[p]}=
a^{[p]}+ ( \chi (a)^p + \alpha^p)c \label{l-chi}
\end{equation}
The original 
construction by Hochschild  \cite{ho1,sul}
uses a $p$-structure \newline 
$(a+ \alpha c)^{[p]}=a^{[p]}+ \chi(a)^pc$. 

\subsubsection{Proposition}
Formula~(\ref{l-chi}) defines a restricted Lie algebra
structure.

\subsubsection{Proof} 
The operation that we define is obviously $p$-linear,
i.e.  $(\beta a+ \beta \alpha c)^{[p]}= 
\beta^p (a+ \alpha c)^{[p]}$.
Let us denote $\mbox{ad}_\fl$ 
by ad and $\mbox{ad}_{\fl_\chi}$
by Ad.
Since $c$ is central, 
$\mbox{Ad}(a+ \alpha c)= \mbox{Ad}\, a$.
Thus, 
$$\mbox{Ad}(a+ \alpha c)^{[p]}= \mbox{Ad}\, a^{[p]} =
\left( \begin{array}{cr}
\mbox{ad}\, a^{[p]} & 0 \\
0 & 0
\end{array} \right )=
\left( \begin{array}{cr}
 (\mbox{ad}\, a)^p & 0 \\
0 & 0
\end{array} \right )=
\left( \begin{array}{cr}
\mbox{ad}\, a & 0 \\
0 & 0
\end{array} \right )^p=$$
$$= ( \mbox{Ad}\, a)^p= (\mbox{Ad} (a + \alpha c))^p.$$
Introducing 
an independent variable $T$, we set $ns_n (a,b)$
to be a coefficient at 
$T^{n-1}$ of $(\mbox{ad}(aT+b))^{p-1}(a)$.
By $S_n$ we denote 
the result of the similar procedure performed
in $\fl_\chi$. It is clear that 
$S_n(a+ \alpha c ,b+\beta c)= s_n (a,b)$.
Finally, $((a+ \alpha c)+ (b + \beta c))^{[p]}=
(a+b)^{[p]}+ 
( \chi (a)^p + \alpha^p+ \chi (b)^p + \beta^p)c=
a^{[p]} +b^{[p]}+ \sum_{i=1}^{p-1}  s_i (a,b)+ 
( \chi (a)^p + \alpha^p+ \chi (b)^p + \beta^p)c=
(a+ \alpha c)^{[p]} + (b+ \beta c)^{[p]}+ 
\sum_{i=1}^{p-1}  S_i (a+\alpha c,b+\beta c).$
$\ \ \ \ \Box$

\subsubsection{When is $\fl_\chi$ split?}
The extension $\fl_{\chi}$ 
is split as an extension of Lie algebras
but not necessarily
as an extension of restricted Lie algebras.

\subsubsection{Lemma} 
The splittings of the extension $\fl_{\chi}$ 
are in one-to-one correspondence 
with $\beta \in \fl^\ast$
satisfying the equations
\begin{eqnarray}
\beta ([x,y])=0, \ \label{split}
\\
\beta (y^{[p]})= \chi (y)^p + \beta (y)^p \label{split2}
\end{eqnarray}
for each $x,y \in \fl$.

\subsubsection{Proof}
Any splitting $\fl \rightarrow \fl_\chi$ 
must be of the form
$y \mapsto y + \beta (y)c$ for some 
$\beta \in \fl^\ast$. 
It is a map of Lie algebras 
if and only if $\beta ([\fl,\fl])=0$.
The splitting
preserves the restricted structure 
if and only if  
equation~(\ref{split2}) holds.
$\Box$

\subsubsection{Corollary} \label{iff}
The canonical map 
$\fl\hookrightarrow\fl_\chi$ 
is a restricted Lie algebra splitting 
if and only if $\chi=0$.

\subsubsection{Corollary}
If $\fl$ is perfect (i.e. $[\fl,\fl]=\fl$) 
then $\fl_\chi$ is split if and
only if $\chi =0$.

\subsection{Connection with universal enveloping algebra}

Our next objective is 
to give another, more intrinsic, description
of $\fl_\chi$. 
Let $\psi:
\fl^\ast\rightarrow (\mbox{Spec}\,O)(\KK)$
be the natural map defined by
$(x^p-x^{[p]})(\psi(\chi))=\chi (x)^p$ for all
$x\in\fl$, $\chi\in\fl^\ast$.

\subsubsection{ 
The restricted enveloping algebra of $\g_m$}
The algebra $u(\g_m)$ is semisimple \cite{ho1}.
Let $c$ be a basis element of $\g_m$
such that $c^{[p]}=c$. 
The elements $1,c, \ldots , c^{p-1}$
form a basis of $u(\g_m)$. 
Let us define {\em the Nielsen polynomial}
$\Ni$ for $\eta \in {\Bbb F}$:
$$\Ni=  \left\{ \begin{array}{ll}
- \sum_{n=1}^{p-1} \frac{c^n}{\eta^n} 
& \mbox{if } \eta \neq 0 \\
1-c^{p-1} & \mbox{if } \eta =0
\end{array} \right. $$
The elements 
$\Ni$ form a complete system of orthogonal idempotents
of $u(\g_m)$.
The idempotent $\Ni$ 
corresponds to the character $\r$,
defined by 
$\r (c )=\eta$, 
of $G_m^1$  
and $\r (\Ni)=1$.

\subsubsection{Theorem} Let $\chi$ 
be a non-zero element of $\fl^\ast$ 
and $\nu=\psi(\chi)$.
Then $U(\fl)\otimes_O\OO({\Bbb F}\nu)$
is isomorphic to $u(\fl_{\chi})$ as a Hopf algebra.

\subsubsection{Proof} 
The map of Lie algebras $\fl \longrightarrow \fl_\chi$
given by $a \mapsto a+0c$ can 
be extended to a map of Hopf algebras
$\zeta :U(\fl) \longrightarrow u(\fl_\chi )$. 
Since $\chi \neq 0$ the algebra 
$u(\fl_\chi )$ is generated by $\fl$
and the map $\zeta$ is onto. 

 On the other hand, there is a natural 
surjective linear map 
$$
\theta :U(\fl) \longrightarrow U(\fl) \otimes_{O} 
{\mathcal O}({\Bbb F} \nu)
$$
given by $y \mapsto y \otimes 1$. 
It follows from Proposition
\ref{hopf} that
$\theta$ is a Hopf algebra map. 
The kernel of $\theta$ is
generated by some elements of $O$. 
It suffices to show
that for each $x \in O$ 
such that $\theta (x)=0$ it holds that
$\zeta (x)=0$. Indeed, 
this condition will imply that $\mbox{ker} \,
\theta \subseteq \mbox{ker} \, \zeta$. 
Thus, there exists a Hopf
algebra map $\kappa : U(\fl) \otimes_{O} 
{\mathcal O}({\Bbb F} \nu) \rightarrow u(\fl_\chi )$. 
It is surjective  since
so is $\zeta$. But both algebras 
have the same dimension $p^{N+1}$ where
$N$ is the dimension of $\fl$. 
Thus, $\kappa$ is an isomorphism.

Let $l_i$ be a basis of $\fl$. Then any $x \in O$ 
has a unique representation
as a polynomial in $l_i^p-l_i^{[p]}$. 
We assume $x= \sum_{\mathfrak k} a_{\mathfrak k}
(l^p-l^{[p]})^{\mathfrak k} \in \mbox{ker} \, \theta$ 
in multi-index notation. This means that
$\sum_{\mathfrak k} a_{\mathfrak k} 
\eta^{\mid {\mathfrak k} \mid} 
(\chi (l)^p)^{\mathfrak k} =0$
for each $\eta \in {\Bbb F}$. 
We can notice that $\zeta (x)=
\sum_{\mathfrak k} a_{\mathfrak k} 
c^{\mid {\mathfrak k} \mid} 
(\chi (l)^p)^{\mathfrak k} \in u(\fl_\chi)$.
Finally, $u(\g_m)$ is a central 
semisimple subalgebra of $u(\fl_\chi )$.
The Pierce decomposition of $u(\fl_\chi )$ is 
$u(\fl_\chi )= \oplus_{\eta \in {\Bbb F}} u(\fl_\chi) \Ni$.
Thus, $\zeta (x) \Ni = 
\sum_{\mathfrak k} a_{\mathfrak k} 
\eta^{\mid {\mathfrak k} \mid} 
(\chi (x)^p)^{\mathfrak k} =0$
for every $\eta \in {\Bbb F}$ and, 
therefore, $\zeta (x)=0$. $\ \ \ \Box$

\subsubsection{}
The theorem clearly fails for $\chi =0$. However,
if one thinks that ${\Bbb F}\cdot 0$ is not just a point
but some infinitesimal neighborhood then the theorem
is adjustable to the case of $\chi =0$. For instance,
the next corollary holds for every $\chi$.

\subsubsection{Corollary} $u(\fl_{\chi})$
is isomorphic to $\oplus_{i \in {\Bbb F}} U_{i \chi}(\fl)$
as an algebra.

\subsubsection{}
A representation of $\fl_\chi$
has {\em a  central charge} $\eta\in{\Bbb F}$  
if
$c$  acts by $\eta$. 
Representations of $\fl$ affording $\chi$ are in one-to-one
correspondence with restricted representations of $\fl_{\chi}$
with a central charge 1.

The next corollary provides an intrinsic construction
of $\fl_\chi$. Recall that the set of primitive elements
of a Hopf algebra $H$ is $P(H)= \{h \in H \mid
\Delta{h}= 1 \otimes h+h \otimes 1 \}$. 
The corollary follows
from the fact that $P(u(\fl))=\fl$.

\subsubsection{Corollary}
$\fl_\chi \cong 
P(U(\fl) \otimes_O {\mathcal O}({\Bbb F} \nu))$. 

\subsubsection{}
One can describe properties of $\fl_\chi$ 
starting from the construction of $\fl_\chi$
as the set of primitive elements of
the Hopf algebra $ U(\fl) \otimes_O 
{\mathcal O}({\Bbb F} \nu)$.
The natural map 
$U(\fl) \otimes_O {\mathcal O}({\Bbb F} \nu) 
\rightarrow u(\fl)$, restricted 
to the set of primitive elements, is the extension
map $\fl_\chi \rightarrow \fl$. 
This extension 
has a canonical Lie algebra splitting 
that  does not
preserve the restricted structure: 
$\fl \hookrightarrow U(\fl ) \rightarrow
U(\fl) \otimes_{O(\fl)} {\mathcal O}({\Bbb F} \nu)$ 
has an image in $\fl_\chi$.
The element $c$ is also canonical: it is easy 
to see that for each $x \in \fl$ such that
$\chi (x) \neq 0$, 
the element 
$c=\frac{x^p-x^{[p]}}{\chi (x)^p} \in U(\fl) 
\otimes_O {\mathcal O}({\Bbb F} \nu)$ 
is central and independent of $x$. It  belongs 
to $\fl_\chi$ but not to  $\fl$.

\subsection{Harish-Chandra pairs}

\subsubsection{} \label{nogroup}
A natural question is to try to find
a central extension 
of algebraic groups 
$G_m \rightarrow \hat{L}_\chi
\rightarrow L$
affording $\fl_{\chi}$ on the tangent level.
There is no such central extension
for a non-zero $\chi$ and a semisimple group
$G$
because all central extensions of $G$ are finite.

\subsubsection{} \label{h-c}
For a nilpotent $\chi\in\g^\ast$,
it is possible to add a piece
of an algebraic group
obtaining a restricted Harish-Chandra pair.
Let us consider an algebraic group $S=St_G( \chi )$,
the stabilizer of $\chi$ in $G$.
The centralizer $C_\g(\chi)$ of $\chi$
contains 
the Lie algebra ${\frak s}$ of $S$. 
We define an embedding of Lie algebras 
$\theta:{\frak s}\hookrightarrow \g_\chi$ 
through the chain
of embeddings 
\begin{equation}
{\mathfrak s}\hookrightarrow C_\g(\chi) \hookrightarrow
\g \hookrightarrow \g_\chi.
\end{equation}
Using the left adjoint action of $G$,
we define an action of $S$ on $\g_\chi$ by
\begin{equation}
g\cdot (x\oplus \alpha c) = (g\cdot x) \oplus \alpha c
\end{equation}
for any $g\in S$, $x\in\g$, $\alpha\in\KK$.
We have to check the following three items
to prove that it is a restricted Harish-Chandra pair.

1. {\em The embedding $\theta$ is of restricted
Lie algebras.}
We need the assumption that
$\chi$ is nilpotent,  
 which means
that $C_{\g}( \chi ) \subseteq \mbox{ker} ( \chi )$
by the definition of a nilpotent element.
Given $a \in C_{\g}( \chi )$, we observe that
$(a\oplus 0c)^{[p]}=
a^{[p]} \oplus \chi (a)^p c = a^{[p]}$.

2. {\em $S$ acts on $\g_\chi$ 
by restricted Lie algebra automorphisms.}
Given $a \in C_{\g}( \chi )$, $g\in S$, 
and $\alpha\in\KK$, 
we observe that 
$(g\cdot(a\oplus \alpha c))^{[p]} =
(g\cdot a)^{[p]} \oplus (\chi (g\cdot a)^p+\alpha^p)c = 
g\cdot (a^{[p]} \oplus ((g^{-1}\cdot \chi)(a)+\alpha^p)c) = 
g\cdot (a+\alpha c)^{[p]}$.

3. {\em The actions of ${\frak s}$ on $\g_\chi$,
induced by the action of  $S$ and the embedding
${\frak s}\hookrightarrow \g_\chi$, 
are the same.}
It is true because the  representation of
$\g$ on $\g_\chi$ is the sum of trivial and adjoint
representations.

\section{Frobenius morphism}

The main object of study in this section is 
a Noetherian algebraic scheme
$X$ over $\KK$. We view $X$ from the two viewpoints.
On the one hand, $X$ is a ringed topological space.
On the other hand, $X$ is a functor
mapping a commutative $\KK$-algebra $R$ to
the set $X(R)$ of points over $R$.

\subsection{Properties of Frobenius morphisms}

\subsubsection{Definition}
Let $X^{(n)}$ be $X$ as 
a scheme (i.e. $X^{(n)}=X$ as a topological space
and $\OO_X^{(n)}=\OO_X$
as a sheaf of rings)
with the new structure over the field:
$
X \longrightarrow \mbox{Spec}\,\KK 
\stackrel{{x^{p^{-n}}}}{\longrightarrow}
\mbox{Spec}\,\KK.
$
{\em The Frobenius morphism} $F_X$, defined on the level
of functions by $f\mapsto f^p$, 
is a morphism of $\KK$-schema:
$
F_X: 
X
\longrightarrow
X^{(1)}. 
$

\subsubsection{Frobenius morphism for a smooth scheme}
The Frobenius morphism $F_X$ is never smooth. 
It is flat if and only if $X$ is smooth 
by Kunz theorem \cite{knu}.  
The following proposition is a technical 
fact about the Frobenius morphism, crucial
for  further study. It would be interesting to know
whether Proposition \ref{smooth} 
holds true for some singular
variety. 

Intuitively,  
the proposition says that the Frobenius map
is locally surjective on points over rings. It holds
if one replaces
the Frobenius map by any
faithfully flat finitely presented map. 

\subsubsection{Proposition} \label{smooth} 
Let $X$ be a smooth algebraic variety
and $R$ be a commutative $\KK$-algebra. 
For each $h \in X^{(1)}(R )$ there
exist a faithfully flat finitely presented 
$R$-algebra $\tilde{R}$ and $y \in X(\tilde{R})$ such
that $F_X(\tilde{R})(y)=X^{(1)}( \varphi )(h)$ 
where $\varphi :R \rightarrow \tilde{R}$
is the natural map.

\subsubsection{Proof}
The Frobenius morphism $F_X$
is flat by
the Kunz theorem. 
It is faithfully flat because
it is surjective on the level
of points over $\KK$. 

We assume that 
$X$ is affine without 
loss of generality since the question
is local. Denote $A= {\mathcal O}(X)$; 
the Frobenius morphism is given 
by the $p$-th
power map $F:A^{(1)} \rightarrow A$,
where $A^{(1)}=\OO(X^{(1)})$ by definition. 
The point $h$ is a map 
of $\KK$-algebras 
$A^{(1)} \rightarrow R$. 
The $R$-algebra $\tilde{R}=R \otimes_{A^{(1)}}A$
is faithfully flat by \cite[1.2.2.5]{dem} and 
obviously finitely presented (it has the same generators
and relations over $R$ as $A$ over $A^{(1)})$. 
The natural map $A \rightarrow \tilde{R}$ is the point
$y$ we are looking for. $\Box$

\subsection{Frobenius neighborhoods}
We define Frobenius neighborhoods and consider
Frobenius kernels as an example of this phenomenon.

\subsubsection{Definition}
Let $Y$ be a closed
subscheme of $X$; then $Y^{(1)}$ 
is naturally a subscheme of $X^{(1)}$.
 Our main concern in this section
is the inverse image subscheme 
$F_X^{-1}(Y^{(1)})$, {\em the Frobenius
neighborhood of $Y$ in $X$}. 
We  denote it by $\widehat{Y}$. This notation 
is ambiguous 
because it is  unclear in which $X$ it is taken.

Assume $Y$ is a closed subscheme of an affine scheme
$X$, determined 
by equations $f_1=0, \,\ldots\, , 
f_m=0$. 
The ideal of $\widehat{Y}$ 
is generated by $f_i^p$. Thus, $\widehat{Y}$ lies
in the $p$-th 
infinitesimal neighborhood of $Y$ and contains
the first infinitesimal neighborhood. 

\subsubsection{Frobenius kernels}
An interesting choice of $X$ and $Y$ 
is $X=L$, 
an algebraic group, and $Y= \{ e \}$, the reduced
identity element. The functoriality of Frobenius
morphism implies that $L^{(1)}$ is an algebraic
group and $F_L$ is a map 
of algebraic group schema. 
The neighborhood $\widehat{Y}$
is the kernel of $F_L$, which is an infinitesimal
finite group scheme 
(called {\em the first Frobenius kernel}).
It will be denoted $L^1$. 

Since $\mo(L^1) \cong u(\fl)^{\ast}$ \cite[1.9.6]{jan1}, 
the $u(\fl)$-modules 
coincide with $\mo(L^1)$-comodules,
i.e. with $L^1$-modules. 

\subsubsection{Frobenius neighborhoods 
in an $L$-variety}
The Frobenius kernel $L^1$ acts on the Frobenius
neighborhood $\widehat{Y}$ of any subvariety $Y$ because
of the functoriality of the Frobenius morphism. 
To prove this, 
pick $g\in L^1(R)$ and $x\in \widehat{Y}(R)$.
We have to show that $gx\in\widehat{Y}(R)$.
The latter means that $F_Z(gx)\in Y^{(1)}(R)$.
But $F_Z(gx)=F_L(g)F_Z(x)$ because of
the functoriality. Now we finish the computation
$F_Z(gx)=F_L(g)F_Z(x)=1_LF_Z(x)=F_Z(x)\in Y^{(1)}(R).$

\subsubsection{Central extensions $L^\chi$ of Frobenius kernels}
It is interesting that $U_\chi (\fl)$-modules
can be also understood in a similar spirit.
They are  representations of a certain central
extension of $L^1$.
The central extension
$$0 \rightarrow \g_m \rightarrow \fl_{\chi}
\rightarrow \fl \rightarrow O$$
gives rise to an exact sequence in the category
of Hopf algebras
$${\KK} \rightarrow u(\g_m) 
\stackrel{\alpha}{\longrightarrow} u(\fl_{\chi})
\rightarrow u(\fl) \rightarrow {\KK}.$$
It is central in a sense that $u( \g_m)$
lies in the center of $u(\fl_{\chi})$. 
The kernel of $\alpha$ is 
an ideal generated by $\g_m$ inside
$u(\fl_{\chi})$. 
We dualize the sequence: 
\begin{equation}
{\KK} \leftarrow u( \g_m)^{\ast} 
\stackrel{\,\beta}{\longleftarrow} 
u(\fl_{\chi})^{\ast}
\leftarrow u(\fl)^{\ast} \leftarrow {\KK}.
\label{functions}
\end{equation}
The centrality of $\g_m$ in $\fl_\chi$
amounts to the fact that 
$b_1 \otimes \beta (b_2) = \beta (b_1) \otimes b_2$
for each 
$b \in u(\fl_{\chi})^{\ast}$. 
The algebra extension
$u(\fl_{\chi})^{\ast}
\supseteq u(\fl)^{\ast}$
is $u(\g_m)^\ast$-Galois \cite{mon,rum}. 
Noticing that $u(\g_m)^\ast\cong\KK{\mathbb Z}_p$,
the Galois condition means
that $u(\fl_{\chi})^{\ast}$ 
is a ${\mathbb Z}_p$-graded algebra such that
$ u(\fl_{\chi})^{\ast}_s \,  
u(\fl_{\chi})^{\ast}_{s^{-1}} = 
u(\fl_{\chi})^{\ast}_e = u(\fl)^{\ast}$
for all $s\in {\mathbb Z}_p$
where $e\in {\mathbb Z}_p$ is the identity element
\cite{mon}. 
Applying the functor Spec to sequence~(\ref{functions}),  
we arrive at 
a central extension of finite infinitesimal
group schema:
$$1 \rightarrow G_m^1 \rightarrow L^{\chi} 
\stackrel{\pi}{\longrightarrow} L^{1} \rightarrow 1$$
where $L^\chi$ is the spectrum 
of $u(\fl_{\chi})^{\ast}$, 
by definition.

\subsubsection{Lemma} \label{f}
For each $\eta\in\BF$, there exists an invertible
element $f\in u(\fl_{\chi})^{\ast}$ such that
$f(xa)=a^\eta f(x)$ for all $a\in G_m^1(R)$,
$x\in L^\chi(R)$, and any
 commutative $\KK$-algebra $R$
(note that $a^\eta$ is well-defined since
$a\in G_m^1(R)=\{r\in R\mid r^p=1\}$). 

\subsubsection{Proof}
The element $\rho=\r\in u(\g_m)^{\ast}$ is group-like
(i.e. $\Delta(\rho)=\rho\otimes\rho$) 
since it is a representation.
Rewriting
$f(xa)=a^\eta f(x)=f(x)\rho(a)$,
we realize that we are looking for an invertible element
$f$ such that $f_1\otimes \beta(f_2)=f \otimes \rho$,
i.e. $f$ is homogeneous of degree $\rho$.
The algebra $u(\fl_{\chi})^{\ast}$ is local. 
As a result, $f$ is invertible
if and only if $\varepsilon (f)\neq 0$. 
The Galois condition 
\cite[Theorem 8.1.7]{mon} implies that 
\begin{equation}
1\in u(\fl_{\chi})^{\ast}_\rho \:
u(\fl_{\chi})^{\ast}_{\rho^{-1}}
\label{graded}
\end{equation} 
where $u(\fl_{\chi})^{\ast}_\rho$ denotes the subspace of
$\rho$-homogeneous elements.
If no such $f$ exists then 
$\varepsilon (u(\fl_{\chi})^{\ast}_\rho)=0$,
which contradicts~(\ref{graded}).
$\ \ \Box$

\subsubsection{Harish-Chandra pairs}
The Harish-Chandra pair $(S,\g_\chi)$, 
constructed in \ref{h-c},
is a central extension of
another pair $(S,\g)$, which 
can be interpreted 
as a Frobenius neighborhood $\widehat{S}$
of $S$ in $G$ since they have the same
categories of representations.
Similarly, one can interpret
the pair $(S,\g_\chi)$ as
a central extension of
$\widehat{S}$ 
by $G_m^1$.

\section{Groupoids}

\subsection{Basics}
We discuss groupoids and
their relevance to Frobenius neighborhoods.
We follow \cite{mac}
for groupoid and Lie algebroid terminology.

\subsubsection{Groupoid scheme}

{\em A groupoid $J$ over a scheme $X$}
is a scheme $J$ over $X\times X$, equipped with morphisms 
$$
m:J^{[2]}=J\times_XJ\rightarrow J, \ 
\iota :X\rightarrow J, \ 
^{-1}:J\rightarrow J
$$
of
multiplication, 
identity
that is a closed embedding, 
and inversion 
such that for any commutative ring $R$ 
the set $J(R)$ is a groupoid
with the base $X(R)$ under 
the structure maps $m(R)$, $\iota(R)$, and
$^{-1}(R)$. Moreover, for any
algebra homomorphism $\mu:R\rightarrow R^\prime$,
the map $J(\mu):J(R)\rightarrow J(R^\prime)$ must be a map
of groupoids.
If the $X\times X$-structure
on $J$ is given by 
$({\mathfrak A},{\mathfrak P})
:J\longrightarrow X\times X$ then the fiber product
$J^{[2]}=J\times_XJ$ is taken using $\fP$ in 
the first position and $\fA$ in the second position.

\subsubsection{Quotients}
A groupoid $J$ over $X$ acts on an $X$-scheme $Y$
if a morphism 
$$\star :(J\stackrel{\fP}{\longrightarrow} X)\times_X Y\rightarrow Y$$
is given 
satisfying associativity and unitarity conditions.
For any
$\KK$-algebra $R$,
an equivalence relation $\sim$ on $Y(R)$ 
is 
$$x\sim y \Longleftrightarrow \exists g\in J(R)\mid 
g\star x=y.
$$
Then $Y/J$ is a sheaf  in
the flat topology on the category of $\KK$-algebras 
associated to the presheaf
$
R\mapsto Y(R)/\sim. 
$

If $Y=X$ and $\star= \fA$ then 
$X/J$ is  a quotient by 
a groupoid as defined in \cite{dem}.

\subsubsection{Action groupoid} \label{action} 
A group scheme $L$ 
action on a scheme $Y$
gives rise to 
{\em the action groupoid} $J_X$
for each closed subscheme $X$ of $Y$.
Note that $X$ need not 
be invariant under the $L$-action.
If $a:L\times Y\rightarrow Y$ is the action map then
$J_X$  is the inverse image scheme:
$J_X=(a|_{L\times X})^{-1}(X)$.
In other words,
${J_X}(R )= \{(g,x) \in L(R) \times X(R ) 
\mid g\cdot x\in X(R)\}$.
The product $m((g,x),(h,y))=(gh,y)$ 
is defined whenever $x=h\cdot y$.

\subsubsection{Product groupoid}
Given a groupoid $J$ over a scheme $Y$ 
and a group scheme $L$,
one can form {\em a product groupoid} $J\times L$ over $Y$.
It is the product scheme 
with the structure maps
$$
\fA^\prime(g,l)=\fA(g),\ \fP^\prime(g,l)=\fP(g),\ 
m^\prime((g,l),(g^\prime,l^\prime))=
(m(g,g^\prime),ll^\prime),\ 
$$
$$
\iota^\prime(x)=(\iota(x),1_L),\ (x,l)^{-1}=
(x^{-1},l^{-1})
$$
for all $g,g^\prime\in J(R)$, $l,l^\prime\in L(R)$, and
$x\in Y(R)$.

\subsubsection{Central extension of a groupoid} \label{gr-ext}
{\em A central extension by 
an Abelian group scheme 
$A$ of groupoid $J$ over $X$} is
a quotient
map  $\pi:J^\prime \rightarrow J$ 
that is a morphism of groupoids
over $\ii_X$. Moreover, an isomorphism must be given 
between the kernel $\pi^{-1}(\iota (X))$ 
and the group scheme $A\times X$, and
the following centrality condition holds.
The equality
\begin{equation} \label{central} 
m(g, (a,\fP(g)))=m((a,\fA(g)),g)
\end{equation} 
must hold
for each $g\in\LL^\chi(R)$, $a\in A(R)$.

\subsubsection{Example}
Let an algebraic group $L$ act on an algebraic variety $Y$.
The central extension $L^\chi$ 
acts on $Y$ 
through $L^\chi\rightarrow L^1\hookrightarrow L$.
For each $X$, a subscheme of $Y$,  
the action groupoid $\LL_X^\chi$
of $L^\chi$ is a central extension 
of the action groupoid $\LL_X^1$ of $L^1$:
\begin{equation} \label{c-gr}
G_m^1\times X \rightarrow \LL_X^\chi \rightarrow 
L_X^1.
\end{equation}

\subsubsection{Proposition} \label{quotient}
If $Y$ is a homogeneous $L$-variety and $X$ is
a smooth subvariety then
$\widehat{X}$ is isomorphic 
to both the quotient of $L^1 \times X$
by the groupoid $\LL_X$ and 
the quotient of $L^\chi\times X$
by the groupoid $\LL_X^\chi$ 
for each $\chi\in\fl^\ast$.

\subsubsection{Proof} 
There is at least one point in $Y(\KK)$ since 
$\KK$ is algebraically closed. 
Thus, we can assume that $Y=L/H$ for some closed subgroup
$H$.  
To treat 
$L^\chi$ and $L^1$ together,
we speak about an infinitesimal group scheme $L^?$
and a groupoid $\LL_X^?$.  
First we show that the action is a quotient map and
then we write down an action of the groupoid $\LL_X^?$.

Thinking of schemes as functors from 
the category of $\KK$-algebras
to the category of sets, we notice that 
the image of the action $L^?\cdot X$ 
is a subfunctor of $\widehat{X}$.
We need to show that it 
is a ``plump'' subfunctor \cite[3.1.1.4 ]{dem},
which means that $\widehat{X}$ 
is a sheaf associated to $L^?X$.
Reiterating the argument before 
Proposition \ref{quotient},
we notice that 
$F_Y(gx)=F_L(\overline{g})F_Y(x)$
for each $g \in L^?(R )$ (where $\overline{g}$
is the image of $g$ in $L^1(R)$), $x \in X(R )$ 
and every $\KK$-algebra $R$. 
Thus, $\widehat{X}(R ) \supseteq L^?(R )\cdot X(R )$.
If $\widetilde{L^?X}$ 
is a sheaf associated to $L^?X$ then 
$\widehat{X} (R )\supseteq 
\widetilde{L^?X} (R )\supseteq L^?(R )X(R )$
for any ring $R$ since $L^?\cdot X$ 
is a subfunctor of a sheaf.

Let us pick $y \in \widehat{X}( R)$. 
We will construct a chain
of faithfully flat finitely presented algebras 
$R\rightarrow R_1\rightarrow 
R_2\rightarrow R_3\rightarrow R_4$ and elements
$g_4\in L^?(R_4)$ and $x_4\in X(R_4)$ 
such that $y_4=g_4x_4$ where
$y_i=\widehat{X}(\pi_i)(y)$ for $\pi_i:R\rightarrow R_i$.
This proves that the action 
$L^?\times X\rightarrow\widehat{X}$
is a quotient map.

By Proposition \ref{smooth}, 
 there exist $R_1$ and $x_1 \in X( R_1)$
such that $F(x_1)=F(y_1)$.
We should notice that this 
is the place that we use the assumption
of $X$ being smooth. 
By the definition of $L/H$, there exist $R_2$
and $a_2,b_2\in L(R_2)$ such that 
$y_2=a_2H(R_2)$ and $x_2=b_2H(R_2)$.
The elements $F_L(a_2)$ and $F_L(b_2)$ 
lie in the same coset. 
Thus, there exists $z_2 \in H^{(1)}(R_2)$
such that $F_L(a_2)=F_L(b_2)z_2$. 
By Proposition \ref{smooth}
used for $H$ 
(any algebraic group is smooth!), there exist 
$R_3$  and $h_3\in H(R_3)$ such that
$z_3 = F_H(h_3)=F_L(h_3)$.
Thus, $F(a_3)=F(b_3)F(h_3)=F(b_3h_3)$.
Let $f_3=a_3h_3^{-1}b_3^{-1}$. It is clear that
$f_3\in L^1(R_3)$ and $f_3x_3=y_3$.

If $?=1$ then we set $R_4=R_3$ and $g_4=f_3$.
If, on the other hand, $?=\chi$ then there exists
$g_4\in G^\chi(R_4)$ such that $\overline{g_4}=f_4$
since $G^\chi\rightarrow G^1$ is a quotient map.
This proves that $\widehat{X}= \widetilde{L^?\cdot X}$.

We define an $\LL^?_X$-action on $L^? \times X$ by
$$(g,x) \star (h,x)= (hg^{-1},gx).$$
$L^?X$ is a quotient functor of $L^? \times X$ by 
the relation 
$$(g,x) \sim (h,y)
\ \ \Longleftrightarrow \ \  gx=hy.$$
 But 
it is equivalent to the condition 
$(h,y) = t \star (g,x)$
where $t=(h^{-1}g,x) \in {\mathcal L}_X^?(R )$. 

Thus, $L^?\cdot X$ is a quotient functor of 
$L^? \times X$ by the groupoid ${\mathcal L}_X$
or $\LL^\chi_X$ correspondently.
This implies that
$\widehat{X}= \widetilde{L^?\cdot X}$ is a quotient sheaf   
$L^? \times X/ {\mathcal L}^?_X$ on the category of $\KK$-algebras.
$\ \ \ \ \ \ \Box$

\subsection{Lie Algebroids}
We discuss Lie theory of groupoids.

\subsubsection{Definition}
Intuitively, 
a Lie algebroid is a tangent structure
to a groupoid \cite{mac}. 
In positive characteristic,
such structure 
is equipped with a $p$-th power map that was axiomatized
by Hochschild \cite{ho3}. 

{\em A restricted Lie algebroid} ${\Bbb L}$ 
on a scheme $X$ is a quasicoherent
$\OO_X$-module that carries a structure 
of a sheaf of restricted Lie algebras
over $\KK$. 
It must be equipped with
an anchor map $\mA:{\Bbb L}\rightarrow TX$ 
that is a morphism of
both $\OO_X$-modules and sheaves of 
restricted Lie algebras. 
Furthermore, it must satisfy the following
identities for  sections 
$u\in\OO_X(V), x,y\in {\Bbb L}(V)$ on an open
subset $V$ of $X$:
\begin{eqnarray}
& & 
[x,uy]=u[x,y]+\mA(x)(u)y, 
\nonumber 
\\ & &  
(ux)^{[p]}=u^px^{[p]}+\mA(ux)^{p-1}(u)x. \label{algebroid}
\end{eqnarray}

For instance, a restricted Lie algebra
is a Lie algebroid over a point.
Another example of a Lie algebroid is the tangent bundle
$TX$. The first relation of (\ref{algebroid}) 
is obvious in this case. The second
one follows from Hochschild's lemma 
\cite[Lemma 1]{ho3}.

\subsubsection{Lie algebroid of a groupoid} 
The Lie algebroid of a groupoid scheme
$J$ over $X$ is the normal sheaf $N_{J|X}$
to the identity morphism $\iota:X\rightarrow J$. 
Quoting \cite{bei}, ``one defines the Lie bracket and projection
by usual formulas'', which one can find in \cite{rum0}. 

\subsubsection{Lie algebroid of an action groupoid}
We consider  a group scheme $L$
acting on a scheme $Y$.
We would like to understand the Lie algebroid $\BL_X$
of the action groupoid of $L$ on $X$ 
for a closed subscheme
$X\subseteq Y$ (see \ref{action}). 
It is easy to see that 
$\BL_Y=\OO_Y\otimes\fl$ as a sheaf with operations
easily computable 
by formulas~(\ref{algebroid}). 

In general, for
on an open affine $V\subseteq X$, 
pick 
an affine open subset 
$V^\prime\subseteq Y$ 
such that $V=X\cap V^\prime$, then  
\begin{equation} \label{tangent}
\BL_X(V) =
\{v|_V \mid
v\in (\OO_Y\otimes\fl) (V^\prime)\, \& \, 
\mA(v) \mbox{ is tangent to }
X\}.
\end{equation}

\subsubsection{Central extensions of Lie algebroids}
Let $L$ be an algebraic group acting on a variety $Y$.
Assume $X$ is a subscheme of $Y$. We have a central
extension~(\ref{c-gr}) of action groupoids 
$\pi :\LL_X^\chi\rightarrow\LL_X$.  
 Their tangent
Lie algebroids ${\BL}_{\chi,X}$ and $\BL_X$ form
a central extension of 
restricted Lie algebroids on $X$:
\begin{equation}
0\rightarrow \g_m\otimes\OO_X\rightarrow 
{\BL}_{\chi,X}
\stackrel{\mbox{\tiny d}\pi }{\longrightarrow}
{\BL}_X \rightarrow 0.
\label{lie}
\end{equation}

\subsubsection{Proposition} \label{bundle}
Let $L$ be a linear algebraic group 
and $Y$ be a homogeneous
$L$-variety. 
For any smooth subvariety $X$,
the Lie algebroid 
$\BL_X$ is a vector subbundle of 
$\OO_X\otimes\fl$.
Similarly, 
$\BL_{\chi,X}$ is a vector subbundle 
of $\OO_X\otimes\fl_\chi$
for  each $\chi\in\fl^\ast$.

\subsubsection{Proof}
To prove the first statement, we show that
the quotient sheaf $\OO_X\otimes\fl / \BL_X$
is locally free. Then $\OO_X\otimes\fl$ is locally
a direct sum of $\BL_X$ and the quotient sheaf since
vector bundles are projective objects in the category
of $\OO$-modules on an affine variety by the
Serre theorem.

Let $J$ be the action groupoid of $L$ on $X$.  
The groupoid $\LL_X$ is the Frobenius
neighborhood of $\iota(X)$ in the groupoid $J$. 
This can be easily
seen because of the functoriality of Frobenius morphism:
points of both $\LL_X$ and $\widehat{\iota (X)}$ 
over $R$ are such 
$(g,x)\in L(R)\times X(R)$ that $F_L(g)=1$. 
Thus, the Lie algebroids
of $J$ and $\LL_X$ coincide, since a normal
bundle is determined by the first order neighborhood that is
contained in the Frobenius neighborhood. 
The quotient sheaf $\OO_X\otimes\fl / \BL_X$
is the normal bundle $N_{L\times X|J}$ restricted to $X$,
which is a subvariety of $J$ under $\iota$.
It suffices to show
that $J$ is smooth, since
a restriction of a locally free sheaf is locally free
and a normal sheaf of an embedding of smooth varieties
is locally free.

Since $Y$ is an $L$-homogeneous variety,
the action morphism 
$a:L\times X\rightarrow Y$ is a submersion and,
therefore, smooth by \cite[Proposition 3.10.4]{har}.
The morphism 
$\fA:J=a^{-1}(X)\rightarrow X$ is smooth, being a base
change of $a$ \cite[Proposition 3.10.1]{har}.
Since $X$ is smooth, then so is $J$.

Now we  prove the second statement.
The sheaf $\BL_{\chi,X}$ is a direct sum of $\BL_X$
and the trivial sheaf $\OO_X\otimes \g_m$. 
Thus, the quotient sheaves
$(\OO_X\otimes\fl) / \BL_X$ and
$(\OO_X\otimes\fl_\chi) / \BL_{\chi,X}$
are the same.
$\ \ \ \Box$

\subsection{Split extensions}

\subsubsection{Definition}
A central  extension of groupoids $\GG^\prime\rightarrow\GG$ 
by an Abelian group scheme $A$
is called {\em split}
if it is isomorphic to the extension 
$\GG\times A\rightarrow\GG$. 

\subsubsection{Lemma \cite{rum0}} \label{split4}
The following statements about a central extension
of groupoids
$\GG^\prime \stackrel{\pi}{\longrightarrow}\GG$
by $A$ 
over a scheme $Y$ are equivalent.
\begin{enumerate}
\item The extension is split. 
\item There exists a groupoid map 
$\mu:\GG\rightarrow\GG^\prime$
such that $\pi\circ\mu=\ii_\GG$.
\item There exists a groupoid map 
$\nu:\GG^\prime\rightarrow A\times X\times X$
such that $\nu(g,x)=(g,x,x)$
for each $(g,x)\in\mbox{ker}\,\pi (R)$.
\item There exists a groupoid map 
$\xi:\GG^\prime\rightarrow A$,
lying  over
the morphism $Y\rightarrow \mbox{Spec}\,\KK$, 
such that $\xi(g,x)=g$
for each $(g,x)\in\mbox{ker}\,\pi (R)$.
\end{enumerate}

\subsubsection{Theorem}
\label{con}
Let a linear algebraic group 
$L$ act on a smooth algebraic variety $Y$
over $\KK$. Let $X$ be 
a smooth subvariety of $Y$ and $\chi\in\fl^\ast$
such that the canonical splitting
of  morphism d$\pi_X$ in (\ref{lie})
is a map of restricted Lie algebras.
Then the central extension~(\ref{c-gr}) 
of action groupoids 
$\pi_X:\fL_X^\chi\rightarrow\fL_X$ is split.

\subsubsection{Proof} A Hopf algebroid $H(J)$ 
of a groupoid $J$ over $X$ is the push-forward
sheaf $(\fA,\fP)_\circ (\OO_J)$.
It is a sheaf of commutative algebras
on $X\times X$, whose local structure
is described in \cite{rum2}.
If the morphism  $(\fA,\fP)$ is affine,
which is the case with action groupoids
of affine group schema,
then the groupoid can be recovered
from its Hopf algebroid as a relative spectrum.
The morphism $\pi_X$ determines a morphism
of Hopf algebroids 
$\pi_X^\#: H(\fL_X)\rightarrow H(\fL_X^\chi)$.
Thus, to split $\pi_X$, it suffices to construct
a morphism of Hopf algebroids splitting $\pi_X^\#$.

The splitting 
of restricted Lie algebroids
determines a morphism of
restricted enveloping $\OO_X$-algebras 
$\zeta :u(\BL_X)\rightarrow u(\BL_{\chi ,X})$.
The left dual morphism $\,^\ast\zeta$ is the splitting
of $\pi_X^\#$ \cite[Corollary 12]{rum2}
because of canonical isomorphisms
$H(\fL_X)\cong\,^\ast u(\BL_X)$
and
$H(\fL_X^\chi)\cong\,^\ast u(\BL_{\chi ,X})$.

The argument in \cite{rum2} is local
but the canonical isomorphisms are 
defined globally
since the construction 
behaves well under localizations.
The ``$O$-good'' condition,
used in \cite{rum2}, 
is that the quotient sheaves
$(\OO_X\otimes\fl) / \BL_X$
and
$(\OO_X\otimes\fl_\chi) / \BL_{\chi,X}$
are locally free.
It is shown in the proof  
of Proposition \ref{bundle}. $\ \Box$

\subsubsection{}
Now we 
choose a connected simply-connected semisimple
algebraic group $G$ as the algebraic group $L$.
The functional $\chi$ is nilpotent. 
The $G$-homogeneous variety is the flag variety $\B$.
$X$ is a subscheme of $\B$. 
We use notation $\pi_X:\GG^\chi_X\rightarrow\GG_X$
for the central extension~(\ref{c-gr}) 
of action groupoids and
d$\pi_X: \BG_{\chi,X}\rightarrow\BG_X$ for the central
extension~(\ref{lie}) of Lie algebroids.

\subsubsection{Infinitesimal splitting condition}
The Lie algebroid
$
{\BG}_{\chi,X}
$
is equal to
$
\BG_X\oplus (\g_m\otimes\OO_X).
$
The inclusion $\gamma_X$ of ${\BG}_X$
into ${\BG}_{\chi,X}$
is a splitting on the level of
Lie algebroids.
We say that 
{\em the infinitesimal splitting condition} 
holds for a subvariety $X$
if $\gamma_X$ is a morphism of restricted Lie algebroids.
The infinitesimal splitting condition
implies that $X$ is 
a subscheme of $\spr$,
which is equivalent to
$\gamma_X$ being a splitting on the diagonal
by Corollary \ref{iff}.
We are going to use the action map $\mu:\g\rightarrow T\B$
in the next proposition.

\subsubsection{Proposition} \label{test}
Let $X$ be a subscheme of $\spr$ such that
the following condition holds
for each Borel subalgebra $\fb\in X(\KK)$:
if $y$ is an element of $\g$ such that the tangent vector 
$\mu(y)_\fb$ defined by $y$ at the point $\fb$
is tangent to $X$ 
then  $\chi(y)=0$.
Under this condition the map
$\gamma_X:\BG_X\rightarrow \BG_{\chi,X}$ is a morphism
of restricted Lie algebroids.

\subsubsection{Proof}
Let $V$ be an open subset of $X$. 
Pick $\sum_iF_i\otimes x_i$ with $F_i\in\OO_X(V)$, $x_i\in\g$
such that $\mA(\sum_iF_i\otimes x_i)$ is tangent to $X$.
Denoting the $p$-th power in $\BG_{\chi,X}$ by $^{(p)}$, we compute
by formulas (\ref{algebroid}). 
\begin{eqnarray*}
& & 
 (\sum_iF_i\otimes x_i)^{(p)}
=
\sum_i (F_i\otimes x_i)^{(p)} + \ldots
=
\\ & & 
\sum_i (F_i^p\otimes x_i^{(p)} + 
\mA(F_ix_i)^{p-1}(F_i)x_i) +\ldots
=
\\  & & 
\sum_i (F_i^p\otimes\chi(x_i)^pc + 
F_i^p\otimes x_i^{[p]} + \mA(F_ix_i)^{p-1}(F_i)x_i) +\ldots
=
\\ & & 
(\sum_i F_i\otimes x_i)^{[p]} + (\sum_iF_i\chi(x_i))^p\otimes c
\end{eqnarray*}
where $\ldots$ denote the terms 
coming from 
the formula for $p$-th degree of a sum
in an associative algebra. 
These terms depend  on 
the adjoint representation only and,
therefore, are the same for $^{(p)}$ and $^{[p]}$.

This argument shows that we have to check that 
$\sum_i F_i\chi(x_i)=0$. 
We check this condition pointwise.
Pick $\fb\in X(\KK)$. Let $y=\sum_i F_i(\fb)x_i\in\g$.
It suffices to deduce $\chi(y)=0$ 
from $\mA(y)$ being tangent to $X$,
which is  
the assumption of the proposition. $\Box$

\subsubsection{Lemma} \label{flag}
Any partial flag variety 
lying in $\spr$ satisfies the infinitesimal splitting condition.

\subsubsection{Proof}
Pick a parabolic subalgebra and 
a Borel subalgebra ${\frak p}
=\mbox{Lie}\, P \supseteq \fb$.
Assume that the partial 
flag variety $X=P\cdot \fb$ lies in $\spr$.
This implies that $\chi$ vanishes on $\frak p$.
But the vector field $\mu(y)$ is tangent to $X$ 
if and only if $y\in {\frak p}$.
We are done by Proposition~\ref{test}.
$\Box$

\subsubsection{}
If $X$ 
is not a partial flag variety  then 
the tangency to $X$ condition is difficult
to put hands on. But if $\mu(y)$ is 
tangent to $X\subseteq \spr$ then it is
also tangent to $\spr$, which implies  
that $\chi([y,\fb])=0$ for each
$\fb\in X(\KK)$. We investigate 
when the latter condition implies $\chi(y)=0$.
 
\subsubsection{Lemma} \label{test3}
If every $\fb\in X(\KK)$ contains an element $h$ such that
ad$^\ast(h)\chi=\chi$
then $X$ satisfies the infinitesimal splitting condition.

\subsubsection{Proof} 
We just need to note that the  pairing
$\g^\ast\times\g\rightarrow\KK$ is $\g$-invariant.
Since
$\chi([y,\fb]) = 0$,
for the choice of $h$ as explained,
$0 = \chi ([y,h])
= \mbox{ad}^\ast h(\chi) (y) = \chi (y). \ \ \ \Box$



\section{Equivariant sheaves and representations}
We introduce the geometric construction of
$\ug$-modules in this section.

\subsection{Equivariant sheaves}
Sheaves equivariant  for groupoids
provide a proper framework
for constructing $U_\chi (\g)$-modules.

\subsubsection{Definition}
We consider a groupoid $J$ over 
an algebraic scheme $Y$.
We notice that 
a groupoid structure gives rise 
to three maps $t_1,t_2,m:J^{[2]}\rightarrow
J$. The maps $t_1$ and $t_2$
are the projections 
to the first and second
component. 
{\em A $J$-equivariant sheaf} is 
an $\OO$-module ${\mathcal F}$ 
on $Y$ with an 
additional structure, namely, an isomorphism 
$I:\fP^\circ(\FF)\rightarrow \fA^\circ(\FF)$ 
of $\OO$-modules on $J$ such that
$$
I|_{i(Y)}:\FF=\fP^\circ(\FF) |_{i(Y)} 
\rightarrow \fA^\circ(\FF) |_{i(Y)}=\FF
$$
is the identity
map and $t_1^\circ I \circ t_2^\circ I=m^\circ I$.
The inverse images are taken in the category of 
$\OO$-modules. 

\subsubsection{Action on fibers}
A $J$-equivariant structure  
gives rise to the action of $J$ 
on the fibers.
Indeed, for each $g \in J$ one obtains an isomorphism
$I_{g}: 
\FF_{\fP(g)}= (\fP^\circ {\mathcal F})_g \rightarrow 
\FF_{\fA(g)}=(\fA^\circ {\mathcal F})_g$.

\subsubsection{}
An $L$-equivariant bundle  ${\mathcal F}$  
may be utilized to construct a large family of 
$L^1$-modules.
Let $X$ be a subscheme of $Y$. 
Then 
$\Gamma ( \widehat{X},{\mathcal F}|_{\widehat{X}})$ 
carries a structure of an $L^1$-module (and, therefore,
a $u(\fl)$-module).

\subsection{Central charge 
of an equivariant vector bundle}

\subsubsection{Definition}
Let $J^\prime \stackrel{\pi}{\longrightarrow}J$
be a central extension by
$G_m^1$
of groupoids over a scheme $X$.
We say that a $J^\prime$-equivariant vector bundle 
$\FF$ has 
{\em a central charge $\eta\in\BF$} if 
$G_m^1$ acts on $\FF$ by the character $\r$.

\subsubsection{} \label{argument} 
If the extension 
of groupoids $\pi:J^\prime\rightarrow J$ is split, 
we can modify a $J$-equivariant structure
on a bundle $\FF$ into a $J^\prime$-equivariant structure
with any central charge $\eta\in\BF$.
Let $\OO_X^\eta$ be the trivial line bundle
with a $G_m^1$-equivariant structure given
by $\r$. 
Thus,
the tensor product $\OO_X^\eta\otimes\FF$
carries a canonical $J^\prime$-equivariant
structure with central charge $\eta$.

\subsubsection{Theorem} \label{splitG}
Let $Y$ be a homogeneous $L$-variety 
and $\FF$ be an $L^1$-equivariant vector bundle on $Y$.
We consider $\chi\in\fl^\ast$ 
and a smooth subvariety $X$ of $Y$
such that
the central extension~(\ref{c-gr}) of action groupoids 
$
\pi_X:{\LL}_X^\chi \longrightarrow \LL_X
$ 
is split.
Then $\FF|_{\widehat{X}}$ admits an $L^\chi$-equivariant 
structure with any central charge $\mu\in\BF$.

\subsubsection{Proof} 
It suffices to exhibit an $L^\chi$-equivariant structure
with a central charge $\mu$ on 
$\OO_{\widehat{X}}$
since
a tensor product of two equivariant vector bundles
has a natural equivariant structure so that central charges
add. Thus, 
$\FF|_{\widehat{X}}\otimes \OO_{\widehat{X}}
\cong \FF|_{\widehat{X}}$ admits an $L^\chi$-equivariant
structure with a central charge $\mu +0$.

By Proposition \ref{quotient},
$\widehat{X}$ is isomorphic to the quotient
$(L^\chi \times X)/\LL_X^\chi$.
The bundle $\OO_{L^\chi\times X}$ admits 
an $\LL_X^\chi$-equivariant structure, called $I$, 
with a central charge $\mu\in\BF$ 
by the argument in \ref{argument} because
the extension $\pi_X$ is split.

The non-trivial part of the proof
is to comprehend 
the quotient $(L^\chi \times X\times \KK)/ \LL_X^\chi$.
The quotient $(L^\chi \times X\times \KK)/G_m^1$
is the trivial line bundle on $L^1 \times X$
because there exists a $G_m^1$-equivariant global
section $s:L^\chi \times X \rightarrow \KK$, 
defined by
$s(g,x)=f(g)$
where a function $f$ is given by Lemma~\ref{f} with
$\eta=-\mu$.
Finally, we observe that 
$(L^\chi \times X\times \KK)/ \LL_X^\chi = 
((L^\chi \times X\times \KK)/G_m^1)/ \LL_X 
\cong (L^1 \times X\times \KK)/\LL_X
\cong \widehat{X}\times\KK$
is the trivial line bundle on $\widehat{X}$,
which inherits
a $L^\chi$-equivariant structure with a central charge
$\mu$
from a $\OO_{L^\chi\times X}.$
$\ \ \ \Box$
 
\subsection{Construction of representations}

\subsubsection{}
We consider a nilpotent functional $\chi\in\g^\ast$.
We say that a subscheme $X$ of $\B$ 
{\em is $\chi$-nice} if it is a smooth subvariety
and satisfies
the infinitesimal splitting condition for $\chi$. 
Every $\chi$-nice subvariety
is a subvariety of the Springer fiber $\spr$
by Corollary~\ref{iff}. 
 
\subsubsection{Theorem} \label{main}
Let us consider
subschema $Z\subseteq X\subseteq\B$
such that $X$ is $\chi$-nice.
If $\FF_\lambda$ is a $G$-equivariant line bundle on $\B$
then the space of sections
$\Gamma (\widehat{Z},\FF_\lambda)$ has a canonical
structure of a $U_\chi(\g)$-module
(the Frobenius neighborhood of $Z$ is taken
inside $\B$).

\subsubsection{Proof} 
Theorem~\ref{con} and Theorem \ref{splitG}
imply that the line bundle $\FF_\lambda |_{\widehat{Z}}$
has a $G^\chi$-equivariant structure with central charge
1.
Therefore, $\Gamma(\widehat{Z},\FF_\lambda)$ 
is a $u(\g_\chi)$-module with central charge 1,
which is canonically a $\ug$-module.
$\ \Box$

\subsubsection{Examples} We would like to compile
a list of known $\chi$-nice subschema.
Any partial flag subvariety in $\spr$ is $\chi$-nice
by Lemma \ref{flag}.
One can check by
a straightforward calculation
that any nilpotent element of ${\mathfrak s}{\mathfrak l}_5$
satisfies the condition of Lemma \ref{test3}.
Thus, this lemma guarantees 
that all smooth subvarieties of $\spr$
are $\chi$-nice for each nilpotent $\chi$ if $\g$ is of type
$A_1$, $A_2$, $A_3$, $A_4$, or $B_2$.

\subsubsection{Stabilizer action}
If $S_1$ is a subgroup of the stabilizer of $\chi$ in $G$
such that
$Y$ is $S_1$-invariant then
$S_1$ also acts on 
the vector space  $\Gamma(\widehat{Y},\FF_\lambda)$. 
It is plausible
that one can combine the actions of $S_1$
and $\g_\chi$ to obtain a representation
of the Harish-Chandra pair $(S_1,\g_\chi)$,
which is a subpair of $(S,\g_\chi)$ constructed 
in \ref{h-c}.

\section{Concluding remarks}

\subsection{Geometric modules}

\subsubsection{The category of geometric modules}
Though the components 
of $\spr$ need not be $\chi$-nice,
we introduce a standard
category of modules. Consider a category $\mathfrak C$
whose objects are pairs $(Z,\lambda)$
where $Z$ is a subscheme of $\B$  
contained in a $\chi$-nice subscheme
and $\lambda$ is a weight. 
The morphism set 
Hom$_{\mathfrak C}((Z,\lambda),(Z^\prime,\lambda^\prime))$
consists of one element if $Z\supseteq Z^\prime$ 
and $\lambda=\lambda^\prime$
and is empty otherwise.
There is a functor
$(Z,\lambda)\mapsto\Gamma (\widehat{Z},\FF_\lambda)$ 
from $\mathfrak C$ to $\ug$-Mod.
A morphism in $\mathfrak C$ goes to
the restriction morphism of the global sections. 
The Abelian subcategory of $\ug$-Mod, generated by the image
of $\mathfrak C$, will be called 
{\em the category of geometric
modules} and denoted $\ug$-Geom.
A module $M$ is called {\em geometric} 
if it is isomorphic
to an object in $\ug$-Geom. 
A filtration (submodule, subquotient) of a geometric
module $M$ is called {\em geometric} if it exists
on  an object 
of $\ug$-Geom
isomorphic to $M$.  

\subsubsection{Question} \label{conjecture}
Are simple $U_\chi(\g)$-modules geometric?

\subsubsection{Parabolic induction}

We want to identify some of 
the geometric modules with modules
constructed by algebraic methods.
Let $P$ be a parabolic subgroup 
containing a Borel subgroup $B$.
Let $U$ be the unipotent radical 
of the opposite  parabolic. $PU$ is a dense
open subset of $G$, isomorphic to $P \times U$. 
It follows that
$\widehat{P/B} \cong P \times U^1$. 
The condition $P/B\subseteq \BB^\chi$ 
is equivalent to $\chi |_{\frak p}=0$
where $\frak p$ is the Lie algebra of $P$. 
The following proposition makes sense
since $u({\frak p})=
U_\chi({\frak p})\subseteq U_\chi(\g)$.

\subsubsection{Proposition} \label{parabolic}
The $U_\chi(\g)$-module 
$\Gamma ( \widehat{P/B}, \FF_{\lambda})$
is isomorphic to 
\newline
$\mbox{Ind}_{U_\chi({\frak p})}^{U_\chi(\g)} 
(\mbox{Ind}_B^P(\KK_{\, -w_0\bullet\lambda}))^\ast$
where $w_0$ is the longest element
of the Weyl group of the Levi factor of $P$ and
$\bullet$ is the dot action.

\subsubsection{Proof}
The Frobenius neighborhood 
$\Sigma$
of the point $P$ in $G/P$
is isomorphic to $G^\chi/P^\chi$. 
The 
$P/B$-bundle 
$\widehat{P/B} 
\stackrel{\pi}{\rightarrow} \Sigma$
is the restriction 
of the natural one to $\Sigma$.
Thus, 
\begin{eqnarray*}
& &
\Gamma (\widehat{P/B},\FF_\lambda)
= \Gamma (\Sigma , \pi_\circ \FF_\lambda) =
\Gamma (\Sigma , G^\chi \times_{P^\chi}
\mbox{Ind}_B^P \KK_\lambda) =
\\ & & 
\mbox{Ind}_{U_\chi({\frak p})}^{U_\chi(\g)} 
((\mbox{Ind}_B^P \KK_\lambda)^\ast)^\ast \cong
\mbox{Ind}_{U_\chi({\frak p})}^{U_\chi(\g)} 
(\mbox{Ind}_B^P(\KK_{\, -w_0\bullet\lambda}))^\ast.
\ \ \Box
\end{eqnarray*}

\subsubsection{The subregular orbit of ${\mathfrak sl}_3$}
We explicate a geometric reason
for a baby Verma module to have more than
one simple quotient.

Let us look at the subregular 
nilpotent orbit of ${\mathfrak sl}_3$.
Let us assume that $p\neq 3$ to identify
$\g$ and $\g^\ast$. 
Choosing a matrix $A$ with $A_{ij}=0$ except
$A_{13}=1$ as $\chi$,
we take the standard Borel subalgebra $\fb$  
to be the intersection
of the two components $Y_1$, $Y_2$
of $\spr$, which is a Dynkin curve
in this case \cite{hum2}.
Now there are non-zero restriction morphisms
$$
\Gamma (\widehat{\spr}, \FF_\lambda) 
\stackrel{i}{\hookrightarrow}
\Gamma (\widehat{Y}_1, \FF_\lambda) \oplus 
\Gamma (\widehat{Y}_2, \FF_\lambda) \rightarrow 
\Gamma (\widehat{\fb}, \FF_\lambda)
$$
for a  weight $\lambda$ 
inside the lowest dominant alcove.
The direct summands in the middle 
are distinct irreducible $\ug$-modules
by Proposition \ref{parabolic} and \cite{kac3,jan2}.
Therefore, the socle of 
$\Gamma (\widehat{\fb}, \FF_\lambda)$ is not simple.
Thus, baby Verma $\ug$-modules 
$\ver$ with this $\fb$, which are isomorphic
to $\Gamma ( \widehat{\fb}, \FF_{-w_0\cdot\lambda})$,
do not have a unique simple quotient.

Another interesting observation is
that  $\Gamma (\widehat{\spr}, \FF_\lambda)$ 
has no natural $\ug$-module structure since the 
embedding $i$ is not an isomorphism. 

\subsection{Deformations of modules}

\subsubsection{}
If $\BB^\chi\subseteq\BB^\eta$ 
then a geometric $U_\chi(\g)$-module
can have a structure of $U_\eta(\g)$-module. 
By Theorem \ref{main}, it suffices to ensure that
 a $\chi$-nice subscheme $Z$ is $\eta$-nice.
Similarly, a geometric
filtration of a $\ug$-module turns out to be a filtration
by $U_\eta(\g)$-modules 
of the corresponding $U_\eta(\g)$-module.

In the particular case of $\eta=0$, 
every geometric $U_\chi(\g)$-module
admits a structure of a restricted $\g$-module
since any smooth subscheme is $0$-nice.
If Question \ref{conjecture} 
has an affirmative answer  then any simple
$U_\chi(\g)$-module 
has a  structure of $u(\g)$-module 
and the dimension of a simple $\ug$-module 
is a sum of dimensions of some simple $u(\g)$-modules.
The case of ${\mathfrak so}_5$ 
has been worked out 
in \cite{rum0}.

\subsubsection{} Let us consider a family of nilpotent
elements $\chi (t)$ and a smooth subvariety 
$Z\subseteq \B$ such that $\B^{\chi (t)}$ contains
$Z$ for each value of the parameter $t$. 
If one can further ensure that $Z$ is $\chi (t)$-nice
for each $t$, 
then we obtain a family of 
$\g$-module structures on the vector space
$\Gamma (\widehat{Z},\FF_\lambda)$ for
each $\lambda\in\Pi$. The $p$-character of
the action at $t$ is $\chi (t)$.

\subsection{Kac-Weisfeiler Conjecture}

\subsubsection{\bf Question} \label{kac}
Let $X$ be a closed subvariety
of a projective algebraic variety $Y$  
and  $\FF$ be a line bundle on $Y$.
When is it true that the dimension of 
$\Gamma (\widehat{X},\FF)$  is divisible by 
$p^{\mbox{\tiny codim}\, X}$?

\subsubsection{}
Affirmative answers to Questions \ref{kac} 
and \ref{conjecture}
would prove the Kac-Weisfeiler Conjecture 
because the dimension formula \cite[6.7]{hum2}
implies that the codimension of $\spr$ in $\B$ 
is equal to $\frac{1}{2}\mbox{dim}\,G\cdot\chi$.

Conversely, if a component of $\spr$ is $\chi$-nice
then the Kac-Weisfeiler conjecture, being the Premet
theorem now \cite{pre}, 
implies an affirmative answer to Question
\ref{kac} for a component of $\spr$ as $X$
and $Y=\B$.
Thus, Question \ref{kac} may be regarded as a geometric
version of the Kac-Weisfeiler conjecture.


\begin{thebibliography}{9977}
\bibitem{bei} A. Beilinson, A. Bernstein, {\em A proof of Jantzen
conjectures}, {Advances in Soviet Mathematics} 
{\bf 16} (part 1) (1993), 1-49.

\bibitem{dem} M. Demazure, P. Gabriel, {\em Groupes Algebriques},
North-Holland Publishing Company, Amsterdam, 1970.

\bibitem{fri} E. M. Friedlander, B. Parshall, 
{\em Modular representation
theory of Lie algebras}, {Amer. J. Math.} 
{\bf 110} (1988), 1055-1094.

\bibitem{har} R. Hartshorne, {\em Algebraic Geometry},
Grad. Texts in Math. {\bf 52}, Springer Verlag, Berlin, 1977.

\bibitem{ho1} G. Hochschild, 
{\em Representations of restricted Lie
algebras of characteristic $p$}, {Proc. Amer. Math. Soc.} 
{\bf 5} (1954), 603-605.
\bibitem{ho2} G. Hochschild, 
{\em Cohomology of restricted Lie algebras},
{Amer. J. Math.} {\bf 76} (1954), 555-580.
\bibitem{ho3} G. Hochschild, 
{\em Simple Lie algebras with purely inseparable
splitting fields of exponent 1}, 
{Trans. Amer. Math. Soc.}
{\bf 79} (1955), 477-489.

\bibitem{hum0} J. E. Humphreys, 
{\em Modular representations of classical Lie
algebras and semisimple groups}, 
{J. Algebra} {\bf 19} (1971), 51-79.
\bibitem{hum2} J. E. Humphreys,  
{\em Conjugacy classes in semisimple algebraic
groups}, Amer. Math. Soc., Providence, 1995. 
\bibitem{hum4} J. E. Humphreys, 
{\em Modular representations of simple
Lie algebras}, 
{Bull. Amer. Math. Soc. (N.S.)} {\bf 35} (1998),
105-122. 

\bibitem{jan1} J. C. Jantzen, Representations of Algebraic Groups,
Academic Press, Orlando, 1987.
\bibitem{jan2} J. C. Jantzen , {\em Subregular nilpotent 
representations of ${\mathfrak s}{\mathfrak l}_n$
and ${\mathfrak s}{\mathfrak o}_{2n+1}$}, 
{Math. Proc. Cambridge Philos. Soc.} {\bf 126} (1999),
223-257.
\bibitem{jan3} J. C. Jantzen, 
{\em Representations of 
${\mathfrak s}{\mathfrak o}_5$ in prime
characteristic}, University of Aarhus 
preprint series {\bf 13}, July 1997.
\bibitem{jan4} J. C. Jantzen, 
{\em Representations of Lie algebras in prime
characteristic}, University of Aarhus 
preprint series {\bf 1}, January 1998.


\bibitem{kac3} V. Kac, {\em On irreducible representations
of Lie algebras of classical type} [Russian],
{Uspekhi Mat. Nauk} {\bf 27} (1972), 237-238.

\bibitem{knu} E. Kunz, {\em
Characterizations of regular local rings
in characteristic $p$}, 
{Amer. J. Math.} {\bf 41} (1969), 772-784.

\bibitem{lus3} G. Lusztig, {\em
Bases in equivariant $K$-theory},
Represent. Theory {\bf 2} (1998), 298-369. 

\bibitem{mac} K. Mackenzie, 
{\em Lie Groupoids and Lie Algebroids in
Differential Geometry}, 
Cambridge University Press, Cambridge, 1987.

\bibitem{mon} S. Montgomery, 
{\em Hopf algebras and their actions on rings},
CBMS Regional Conf. Ser. in Math. {\bf 82}, Amer. Math. Soc.,
Providence, 1993.

\bibitem{nie} G. Nielsen, {\em 
A determination of the minimal right ideals
in the enveloping algebra 
of a Lie algebra of classical type},
Ph.D. dissertation, Univ. of Wisconsin, 1963.

\bibitem{pre} A. Premet, {\em 
Irreducible representations of Lie algebras
of reductive groups and 
the Kac-Weisfeiler conjecture}, {Invent. Math.}
{\bf 121} (1995), 79-117.

\bibitem{rum} D. Rumynin, {\em 
Hopf-Galois extensions with central
invariants and their geometric properties}, 
{Algebr. Represent. Theory}, {\bf 1} (1998), 353-381.
\bibitem{rum0} D. Rumynin, {\em Modular Lie algebras and their
representations}, 
Ph.D. dissertation, 
Univ. of Massachusetts at Amherst, 1998.
\bibitem{rum2} D. Rumynin, {\em Duality for Hopf algebroids},
Journal of Algebra, to appear.


\bibitem{far} H. Strade, R. Farnsteiner, 
{\em Modular Lie algebras
and their representations}, Marcel Dekker, New York, 1988.

\bibitem{sul} J. B. Sullivan, {\em 
Lie algebra cohomology at irreducible
modules}, {Illinois J. Math.} {\bf 23} (1979), 363-373.

\bibitem{kac1} B. Yu. Weisfeiler and V. G. Kac, 
{\em On irreducible representations
of Lie $p$-algebras}, 
{Funct. Anal. Appl.} {\bf 5} (1971), 111-117.

\end{thebibliography}
\end{document}